УДК 517.984.5
# О РАЗРЕШИМОСТИ СМЕШАННОЙ ЗАДАЧИ ДЛЯ УРАВНЕНИЯ ДРОБНОГО ПОРЯДКА С ЗАПАЗДЫВАЮЩИМ АРГУМЕНТОМ ПО ВРЕМЕНИ И ПСЕВДОДИФФЕРЕНЦИАЛЬНЫМИ ОПЕРАТОРАМИ С НЕЛОКАЛЬНЫМИ КРАЕВЫМИ УСЛОВИЯМИ В КЛАССАХ СОБОЛЕВА

Бабаев М.М[1]

**АННОТАЦИЯ.** В данной работе изучена смешанная задача для уравнения дробного порядка с запаздывающим аргументом по времени и псевдодифференциальными операторами, связанными операторами Лапласа с нелокальными краевыми условиями в классах Соболева. Решения начально-граничной задачи построено в виде суммы ряда по системе собственных функций многомерной спектральной задачи. У спектральной задачи найдены собственные значения и построена соответствующая система собственных функций. Показано, что эта система собственных функций является полной и образует базис Рисса в подпространствах Соболева. На основании полноты системы собственных функций доказана теорема единственности решения задачи. В подпространствах Соболева доказано существование регулярного решения поставленной начально-граничной задачи.

**Ключевые слова:** дифференциальное уравнение в частных производных с запаздывающим аргументом, дробная производная по времени, начально-граничная задача, спектральный метод, собственные значения, собственные значения, собственные функции, полнота, базис Рисса, единственность, существование, ряд.

**1. Постановка задачи.** Как известно, что в физике твердого тела изучаются так называемые фрактальные среды, в частности, явления диффузия в них. В одной из моделей, диффузия в сильно пористой среде описывается уравнением типа уравнения теплопроводности, но с дробной производной по временной координате с запаздывающим аргументом. Многие задачи о колебаниях балок и пластин, которые имеют большое значение в строительной механике, приводят к дифференциальным уравнениям [1, с. 141-143], [2, с. 278-280], [3, гл.3].

Отметим также, что к уравнению колебаний балки приходят во многих задачах при расчёте устойчивости вращающихся валов и изучении вибрации кораблей [4-7].

---

[1] **Бабаев Махкамбек Мадаминович** Национальний Университет Узбекистана, г. Ташкент, Узбекистан, babayevm@mail.ru

В данной работе рассматривается дифференциальное уравнение с дробной производной вида

$$D_{0t}^{\alpha}u(x,t) + Bu(x,t) + Cu(x,t-\tau) = f(x,t), \quad (x,t) \in Q = \Pi \times (0,T), \quad \alpha > 0 \quad (1)$$

с начальными функциями

$$\begin{cases} D_{0t}^{\alpha-k}u(x,t)\big|_{t=+0} = \varphi_k(x), \quad x = (x_1,...,x_j,...,x_N) \in \Pi, \quad k=1,2,...,l-1, \\ D_{0t}^{\alpha-l}u(x,t) = \varphi_l(x,t), \quad (x,t) = (x_1,...,x_j,...,x_N,t) \in \Pi \times (-\tau,0), \quad l = -[-\alpha]. \end{cases} \quad (2)$$

Здесь при $\alpha < 0$ дробная интеграл $D^{\alpha}$ имеет вид

$$D_{at}^{\alpha}u(x,t) = \frac{sign(t-a)}{\Gamma(-\alpha)} \int_a^t \frac{u(x,\tau)d\tau}{|t-\tau|^{\alpha+1}},$$

при $\alpha = 0$, то $D_{at}^{\alpha}u(x,t) = u(x,t)$, а при $l-1 < \alpha \leq l$, $l \in \Box$, дробная производная имеет вид

$$D_{at}^{\alpha}u(x,t) = sign^l(t-a)\frac{d^l}{dt^l}D_{at}^{\alpha-l}u(x,t) = \frac{sign^{l+1}(t-a)}{\Gamma(l-\alpha)}\frac{d^l}{dt^l}\int_a^t \frac{u(x,\tau)\cdot d\tau}{|t-\tau|^{\alpha-l+1}}.$$

Оператор

$$Au(x,t) = -\Delta u(x,t) = -\sum_{j=1}^{N}\frac{\partial^2 u(x,t)}{\partial x_j^2}, \ (x,t) \in Q$$

определяется в гильбертовом пространстве $H = W_2^{s_1,s_2,...,s_N}(\Pi)$ с областью определением

$$\overset{0}{W}_2^{s_1,s_2,...,s_N}(\Pi) = \left\{ u(x,t): \quad \alpha_j \cdot \frac{\partial^{2k_j}u(x_1,...,x_j,...,x_N,t)}{\partial x_j^{2k_j}}\bigg|_{x_j=0} + \right.$$

$$+ \beta_j \cdot \frac{\partial^{2k_j}u(x_1,...,x_j,...,x_N,t)}{\partial x_j^{2k_j}}\bigg|_{x_j=\pi} = 0, \ 1 \leq j \leq p, \ 0 \leq k_j < \frac{2s_j - N}{4}, \ s_j > 2 + \frac{N}{2},$$

$$\beta_j \cdot \frac{\partial^{2k_j+1}u(x_1,...,x_j,...,x_N,t)}{\partial x_j^{2k_j+1}}\bigg|_{x_j=0} + \alpha_j \cdot \frac{\partial^{2k_j+1}u(x_1,...,x_j,...,x_N,t)}{\partial x_j^{2k_j+1}}\bigg|_{x_j=\pi} = 0,$$

$$1 \leq j \leq p, \ 0 \leq k_j < \frac{2s_j - N - 2}{4}, \ s_j > 2 + \frac{N}{2},$$

$$\frac{\partial^{2k_j}u(x_1,...,x_j,...,x_N,t)}{\partial x_j^{2k_j}}\bigg|_{x_j=0} = \frac{\partial^{2k_j}u(x_1,...,x_j,...,x_N,t)}{\partial x_j^{2k_j}}\bigg|_{x_j=\pi},$$

$$p+1 \leq j \leq q, \ 0 \leq k_j < \frac{2s_j - N}{4}, \ s_j > 2 + \frac{N}{2},$$

$$\left.\frac{\partial^{2k_j+1}u(x_1,\ldots,x_j,\ldots,x_N,t)}{\partial x_j^{2k_j+1}}\right|_{x_j=0}=\left.\frac{\partial^{2k_j+1}u(x_1,\ldots,x_j,\ldots,x_N,t)}{\partial x_j^{2k_j+1}}\right|_{x_j=\pi},$$

$$p+1\leq j\leq q,\ k_j<\frac{2s_j-N-2}{4},\ s_j>2+\frac{N}{2},$$

$$\left.\frac{\partial^{2k_j}u(x_1,\ldots,x_j,\ldots,x_N,t)}{\partial x_j^{2k_j}}\right|_{x_j=0}=0,\quad \left.\frac{\partial^{2k_j}u(x_1,\ldots,x_j,\ldots,x_N,t)}{\partial x_j^{2k_j}}\right|_{x_j=\pi}=0,$$

$$q+1\leq j\leq N,\ 0\leq k_j<\frac{2s_j-N}{4},\ s_j>2+\frac{N}{2},$$

$$1\leq p\leq q\leq N,\ u(x,t)\in W_2^{s_1,s_2,\ldots,s_N}(\Pi)\Big\} \quad (3)$$

где $\tau,\ T>0$ – постоянные, $l\in\square$, $(x,t)=(x_1,\ldots,x_j,\ldots,x_N,t)\in\Pi\times(0,T)$, $\Pi=(0,\pi)\times\ldots\times(0,\pi)$, $\alpha_i=const,\ \beta_i=const$, $\alpha_j\neq 0$, $\beta_j\neq 0$, $|\alpha_j|\neq|\beta_j|$ при $1\leq j\leq p$.

Для произвольной кусочно равномерно непрерывной на $(-\infty,+\infty)$ вещественной функции $f(\lambda)$ определим псевдодифференциальный оператор $B=f(A)$ со следующими равенствами (см. [8, с. 381])

$$Bu(x,t)=\sum_{n_1=1}^{\infty}\ldots\sum_{n_N=1}^{\infty}f(\lambda_{n_1,\ldots,n_N})T_{n_1,\ldots,n_N}(t)v_{n_1,\ldots,n_N}(x),\quad (x,t)\in Q,$$

где $T_{n_1,\ldots,n_N}(t)=(u(x,t),v_{n_1,\ldots,n_N}(x))_{W_2^{s_1,s_2,\ldots,s_N}(\Pi)}$ – коэффициент Фурье функции $u(x,t)$, с областью определения $D(B)$, состоящей из тех $u(x,t)\in \overset{0}{W}_2^{s_1,s_2,\ldots,s_N}(\Pi)$, для которых

$$\sum_{n_1=1}^{\infty}\ldots\sum_{n_N=1}^{\infty}\left|f(\lambda_{n_1,\ldots,n_N})\right|^2\left|T_{n_1,\ldots,n_N}(t)\right|^2<+\infty.$$

Аналогично, для произвольной кусочно равномерно непрерывной на $(-\infty,+\infty)$ вещественной функции $f(\lambda)$ определим псевдодифференциальный оператор $C=g(A)$ со следующими равенствами (см. [8, с. 381])

$$Cu(x,t-\tau)=\sum_{n_1=1}^{\infty}\ldots\sum_{n_N=1}^{\infty}f(\lambda_{n_1,\ldots,n_N})T_{n_1,\ldots,n_N}(t-\tau)v_{n_1,\ldots,n_N}(x),\quad (x,\ t-\tau)\in Q,$$

где $T_{n_1,\ldots,n_N}(t-\tau) = (u(x,t-\tau), v_{n_1,\ldots,n_N}(x))_{W_2^{s_1,s_2,\ldots,s_N}(\Pi)}$ – коэффициент Фурье функции $u(x,t-\tau)$, с областью определения $D(C)$, состоящей из тех $u(x,t-\tau) \in \overset{0}{W}_2^{s_1,s_2,\ldots,s_N}(\Pi)$, для которых

$$\sum_{n_1=1}^{\infty} \ldots \sum_{n_N=1}^{\infty} \left|g(\lambda_{n_1,\ldots,n_N})\right|^2 \left|T_{n_1,\ldots,n_N}(t-\tau)\right|^2 < +\infty.$$

Далее, достаточно гладкие функции $f(x,t)$, $\varphi_k(x)$ $k=1,2,\ldots l-1$, $\varphi_l(x,t)$ при каждом $t$ разлагаемые функции по собственным функциям $\{v_{m_1,\ldots,m_N}(x_1,\ldots,x_N)\}_{(m_1,\ldots,m_p)\in\mathbb{Z}^p \times (m_{p+1},\ldots,m_q)\in\mathbb{Z}^{q-p} \times (m_{q+1},\ldots,m_N)\in\mathbb{Z}^{N-q}}$ систему собственных функций спектральной задачи:

$$\Delta v(x) + \mu v(x) = 0, \qquad (4)$$

$$\begin{cases}
\alpha_j \cdot v(x_1,\ldots,x_j,\ldots,x_N)\big|_{x_j=0} + \beta_j \cdot v(x_1,\ldots,x_j,\ldots,x_N)\big|_{x_j=\pi} = 0, \ 1 \le j \le p, \\
\beta_j \cdot \dfrac{\partial^{2k_j+1} v(x_1,\ldots,x_j,\ldots,x_N)}{\partial x_j^{2k_j+1}}\bigg|_{x_j=0} + \alpha_j \cdot \dfrac{\partial^{2k_j+1} v(x_1,\ldots,x_j,\ldots,x_N)}{\partial x_j^{2k_j+1}}\bigg|_{x_j=\pi} = 0, \ 1 \le j \le p, \\
v(x_1,\ldots,x_j,\ldots,x_N)\big|_{x_j=0} = v(x_1,\ldots,x_j,\ldots,x_N)\big|_{x_j=\pi}, \quad p+1 \le j \le q, \\
\dfrac{\partial v(x_1,\ldots,x_j,\ldots,x_N)}{\partial x_j}\bigg|_{x_j=0} = \dfrac{\partial v(x_1,\ldots,x_j,\ldots,x_N)}{\partial x_j}\bigg|_{x_j=\pi}, \quad p+1 \le j \le q, \\
v(x_1,\ldots,x_j,\ldots,x_N)\big|_{x_j=0} = 0, \ v(x_1,\ldots,x_j,\ldots,x_N)\big|_{x_j=\pi} = 0, \quad q+1 \le j \le N, \\
1 \le p \le q \le N
\end{cases} \qquad (5)$$

Здесь $\mu_{m_1,\ldots,m_N} = \sum_{j=1}^{p}(2m_j+\varphi_j)^2 + \sum_{j=p+1}^{q}(2m_j)^2 + \sum_{j=q+1}^{N}(m_j)^2$, а числа $\varphi_j = \dfrac{1}{\pi}\arccos\dfrac{-2\alpha_j\beta_j}{\alpha_j^2+\beta_j^2}$. Соответствующая система собственных функций имеет вид

$$\{v_{m_1,\ldots,m_N}(x_1,\ldots,x_N)\}_{(m_1,\ldots,m_p)\in\mathbb{Z}^p \times (m_{p+1},\ldots,m_q)\in\mathbb{Z}^{q-p} \times (m_{q+1},\ldots,m_N)\in\mathbb{Z}^{N-q}} =$$

$$= \left\{\prod_{j=1}^{p} \sqrt{\dfrac{2}{\pi}} \cdot \dfrac{sign(\beta_j^2-\alpha_j^2) \cdot \alpha_j \sin\lambda_{m_j} x_j + \beta_j \cos\lambda_{m_j} x_j}{\sqrt{\alpha_j^2+\beta_j^2} \cdot \sqrt{1+\left|\lambda_{m_j}\right|^{2s_j}}}\right\}_{(m_1,\ldots,m_p)\in\mathbb{Z}^p} \times$$

$$\times \left\{ \prod_{j=p+1}^{q} \frac{1}{\sqrt{\pi}} \cdot \frac{1}{\sqrt{1+|2m_j|^{2s_j}}} exp(i2m_j x_j) \right\}_{(m_{p+1},\ldots,m_q)\in \mathbb{Z}^{q-p}} \times$$

$$\times \left\{ \prod_{j=q+1}^{N} \sqrt{\frac{2}{\pi}} \cdot \frac{1}{\sqrt{1+|m_j|^{2s_j}}} \sin(m_j x_j) \right\}_{(m_{q+1},\ldots,m_N)\in \mathbb{Z}^{N-q}}.$$

Отметим, что аналогичные начально-граничные задачи были исследованы в работе [9] и [10].

Введем пространство $W_2^s(0,l)$ с нормой $\|f\|^2_{W_2^s(0,l)} = \|f\|^2_{L_2(0,l)} + \|D^s f\|^2_{L_2(0,l)}$, где $s$ произвольное натуральное число, при этом $W_2^0(0,l) = L_2(0,l)$.

Скалярное произведение в пространстве $W_2^{s_1,s_2,\ldots,s_N}(\Pi)$, вводится так:

$$(f(x), g(x))_{W_2^{s_1,s_2,\ldots,s_N}(\Pi)} = (f(x), g(x))_{L_2(\Pi)} + \sum_{j_1=1}^{N} (D_{x_{j_1}}^{s_{j_1}} f(x), D_{x_{j_1}}^{s_{j_1}} g(x))_{L_2(\Pi)} +$$

$$+ \sum_{1\le j_1 < j_2 \le N} (D_{x_{j_1}}^{s_{j_1}} D_{x_{j_2}}^{s_{j_2}} f(x), D_{x_{j_1}}^{s_{j_1}} D_{x_{j_2}}^{s_{j_2}} g(x))_{L_2(\Pi)} + \ldots +$$

$$+ \sum_{1\le j_1 < j_2 < \ldots < j_N \le N} (D_{x_{j_1}}^{s_{j_1}} D_{x_{j_2}}^{s_{j_2}} \ldots D_{x_{j_N}}^{s_{j_N}} f(x), D_{x_{j_1}}^{s_{j_1}} D_{x_{j_2}}^{s_{j_2}} \ldots D_{x_{j_N}}^{s_{j_N}} g(x))_{L_2(\Pi)}.$$

Соответственно норма в пространстве $W_2^{s_1,s_2,\ldots,s_N}(\Pi)$, вводится так:

$$\|f\|^2_{W_2^{s_1,s_2,\ldots,s_N}(\Pi)} = \|f\|^2_{L_2(\Pi)} + \sum_{j_1=1}^{N} \left\|D_{x_{j_1}}^{s_{j_1}} f(x)\right\|_{L_2(\Pi)} +$$

$$\sum_{1\le j_1 < j_2 \le N} \left\|D_{x_{j_1}}^{s_{j_1}} D_{x_{j_2}}^{s_{j_2}} f(x)\right\|_{L_2(\Pi)} + \ldots + \sum_{1\le j_1 < j_2 < \ldots < j_N \le N} \left\|D_{x_{j_1}}^{s_{j_1}} D_{x_{j_2}}^{s_{j_2}} \ldots D_{x_{j_N}}^{s_{j_N}} f(x)\right\|_{L_2(\Pi)}.$$

**2. Полнота системы собственных функций в подпространствах Соболева.** Обозначим через $\overset{0}{W}_2^{s_1,s_2,\ldots,s_N}(\Pi)$ множество всех функций $f(x) \in W_2^{s_1,s_2,\ldots,s_N}(\Pi)$, удовлетворяющих граничным условиям (5).

Справедливо следующая

**Теорема 1.** *Пусть* $\alpha_j \ne 0$, $\beta_j \ne 0$, $|\alpha_j| \ne |\beta_j|$ *действительные числа при каждом* $1 \le j \le p$ *и*

$$\rho = \max_{1\le j \le p} \sqrt{\theta_j^2 + 2(\frac{\theta_j}{\sqrt{2}} + (\varphi_j+1)^{s_j} - 1)^2} \cdot \sigma(s_j) < 1,$$

*где* $\sigma(0) = \frac{1}{\sqrt{2}}$, $\sigma(s_j) = 1$, *при* $s_j > 0$, $\theta_j = \sqrt{2} \cdot \max_{x \in [0,\pi]} |e^{i\varphi_j x} - 1|$, $\lambda_{m_j} = 2m_j + \varphi_j$,

$\varphi_j = \frac{1}{\pi} \arccos \frac{-2\alpha_j \beta_j}{\alpha_j^2 + \beta_j^2}$, $m_j \in \mathbb{Z}$. *Тогда система собственных функций*

$$\{v_{m_1 \ldots m_N}(x_1, \ldots, x_N)\}_{(m_1, \ldots, m_p) \in \mathbb{Z}^p \times (m_{p+1}, \ldots, m_q) \in \mathbb{Z}^{q-p} \times (m_{q+1}, \ldots, m_N) \in \mathbb{Z}^{N-q}} =$$

$$= \left\{ \prod_{j=1}^{p} \sqrt{\frac{2}{\pi}} \cdot \frac{sign(\beta_j^2 - \alpha_j^2) \cdot \alpha_j \sin \lambda_{m_j} x_j + \beta_j \cos \lambda_{m_j} x_j}{\sqrt{\alpha_j^2 + \beta_j^2} \cdot \sqrt{1 + |\lambda_{m_j}|^{2s_j}}} \right\}_{(m_1, \ldots, m_p) \in \mathbb{Z}^p} \times$$

$$\times \left\{ \prod_{j=p+1}^{q} \frac{1}{\sqrt{\pi}} \cdot \frac{1}{\sqrt{1 + |2m_j|^{2s_j}}} exp(i2m_j x_j) \right\}_{(m_{p+1}, \ldots, m_q) \in \mathbb{Z}^{q-p}} \times$$

$$\times \left\{ \prod_{j=q+1}^{N} \sqrt{\frac{2}{\pi}} \cdot \frac{1}{\sqrt{1 + |m_j|^{2s_j}}} \sin(m_j x_j) \right\}_{(m_{q+1}, \ldots, m_N) \in \mathbb{Z}^{N-q}},$$

*спектральной задачи (4)-(5) образует полной ортонормированный системой в классах Соболева* $\overset{0}{W}_{2}^{s_1, s_2, \ldots, s_N}(\Pi)$.

**Теорема 2.** *Пусть* $\alpha_j \neq 0$, $\beta_j \neq 0$, $|\alpha_j| \neq |\beta_j|$ *действительные числа при каждом* $1 \leq j \leq p$ *и*

$$\rho = \max_{1 \leq j \leq p} \sqrt{\theta_j^2 + 2(\frac{\theta_j}{\sqrt{2}} + (\varphi_j + 1)^{s_j} - 1)^2} \cdot \sigma(s_j) < 1,$$

*где* $\sigma(0) = \frac{1}{\sqrt{2}}$, $\sigma(s_j) = 1$, *при* $s_j > 0$, $\theta_j = \sqrt{2} \cdot \max_{x \in [0,\pi]} |e^{i\varphi_j x} - 1|$, $\lambda_{m_j} = 2m_j + \varphi_j$,

$\varphi_j = \frac{1}{\pi} \arccos \frac{-2\alpha_j \beta_j}{\alpha_j^2 + \beta_j^2}$, $m_j \in \mathbb{Z}$. $s_j > k + \frac{N}{2}$, $k \geq 0$, $k \in \mathbb{Z}$. *Тогда ряд Фурье функции* $f(x) \in \overset{0}{W}_{2}^{s_1, s_2, \ldots, s_N}(\Pi) \cap C^k(\Pi)$ *по ортонормированных собственных функций*

$$\{v_{m_1 \ldots m_N}(x_1, \ldots, x_N)\}_{(m_1, \ldots, m_p) \in \mathbb{Z}^p \times (m_{p+1}, \ldots, m_q) \in \mathbb{Z}^{q-p} \times (m_{q+1}, \ldots, m_N) \in \mathbb{Z}^{N-q}} =$$

$$= \left\{ \prod_{j=1}^{p} \sqrt{\frac{2}{\pi}} \cdot \frac{sign(\beta_j^2 - \alpha_j^2) \cdot \alpha_j \sin \lambda_{m_j} x_j + \beta_j \cos \lambda_{m_j} x_j}{\sqrt{\alpha_j^2 + \beta_j^2} \cdot \sqrt{1 + |\lambda_{m_j}|^{2s_j}}} \right\}_{(m_1, \ldots, m_p) \in \mathbb{Z}^p} \times$$

$$\times \left\{ \prod_{j=p+1}^{q} \frac{1}{\sqrt{\pi}} \cdot \frac{1}{\sqrt{1+|2m_j|^{2s_j}}} exp(i2m_j x_j) \right\}_{(m_{p+1},...,m_q) \in \square^{q-p}} \times$$

$$\times \left\{ \prod_{j=q+1}^{N} \sqrt{\frac{2}{\pi}} \cdot \frac{1}{\sqrt{1+|m_j|^{2s_j}}} \sin(m_j x_j) \right\}_{(m_{q+1},...,m_N) \in \square^{N-q}} ,$$

*спектральной задачи (4)-(5) сходится по норме пространства* $C^k(\prod)$ *к функции* $f(x)$.

Доказательство теоремы 1 и 2 можно найти в работе [9].

**3. Существование и единственность решения начально-граничной задачи.** Регулярным решением уравнения (1) в области $Q = \prod \times (0,T)$, $T > 0$, назовем функцию $u(x,t)$ из класса $u(x,t) \in C(\overline{Q})$, $D_{0t}^{\alpha-i} u(x,t) \in C(\overline{Q}), i=1,2,...,p$

$$D_{0t}^{\alpha} u(x,t) \in C(Q), \quad \frac{\partial u(x,t)}{\partial x_j} \in C(\overline{Q}), \quad \frac{\partial u^2(x,t)}{\partial x_j^2} \in C(Q), \quad j=1,...,\square \quad , \quad \text{и}$$

удовлетворяющую уравнению (1) во всех точках $(x,t) \in Q$.

Обозначим через $\overset{0}{W}_2^{s_1,s_2,...,s_N;\theta}(Q)$ множество всех функций $u(x,t) \in W_2^{s_1,s_2,...,s_N;\theta}(Q)$, удовлетворяющих граничным условиям (3).

Функцию $u(x,t)$ назовем регулярным решением задача (1)-(3) в области $Q = \prod \times (0,T)$, если функция $u(x,t)$ регулярные решения уравнения (1) в области $Q = \prod \times (0,T)$ и удовлетворяет начальным функциям и граничным условиям (2) и (3).

Пусть функция $u(x,t) \in W_2^{s_1,s_2,...,s_N;\theta}(Q)$ с показателем $s_1 = s_2 = ... = s_N = = 2 + \frac{N}{2}, \theta = -[-\alpha]$ удовлетворяют уравнению (1) во всех точках $(x,t) \in Q$ и удовлетворяет начальным и граничным условиям (2) и (3). Тогда функция $u(x,t)$ является регулярным решением задача (1)-(3) в области $Q = \prod \times (0,T)$.

Введем функции

$$T_{m_1...m_N}(t) = \int_{\prod} u(y,t) \tilde{\upsilon}_{m_1...m_N}(y) dy , \qquad (6)$$

где

$$\tilde{\upsilon}_{m_1...m_N}(x_1,...,x_N) = = \prod_{j=1}^{p} \sqrt{\frac{2}{\pi}} \cdot \frac{sign(\beta_j^2 - \alpha_j^2) \cdot \alpha_j \sin \lambda_{m_j} x_j + \beta_j \cos \lambda_{m_j} x_j}{\sqrt{\alpha_j^2 + \beta_j^2}} \times$$

$$\times \prod_{j=p+1}^{q} \frac{1}{\sqrt{\pi}} exp(i2m_j x_j) \cdot \prod_{j=q+1}^{N} \sqrt{\frac{2}{\pi}} \cdot \sin(m_j x_j) \qquad (7)$$

при $(m_1,...,m_p) \in \square^p, (m_{p+1},...,m_q) \in \square^{q-p}, (m_{q+1},...,m_N) \in \square^{N-q}$.

В силу (1)-(3) неизвестные функции $T_{m_1...m_N}(t)$ удовлетворяют уравнениям

$$D_{0t}^{\alpha} T_{m_1...m_N}(t) + \mu_{m_1...m_N}(a^2 T_{m_1...m_N}(t) + b^2 T_{m_1...m_N}(t-\tau)) = f_{m_1...m_N}(t) \qquad (8)$$

и начальным функциям

$$\begin{cases} D_{0t}^{\alpha-i} T_{m_1...m_N}(t)\big|_{t=+0} = \varphi_{i,m_1...m_N}, \quad i=1,2,...l-1, \\ D_{0t}^{\alpha-l} u(x,t) = \varphi_l(x,t), \ (x,t) = (x_1,...,x_j,...,x_N,t) \in \prod \times (-\tau, 0) \end{cases} \qquad (9)$$

где

$$f_{m_1...m_N}(t) = \int_\Pi f(y,t)\tilde{\upsilon}_{m_1...m_N}(y)dy, \quad \varphi_{j,m_1...m_N} = \int_\Pi \varphi_j(y)\tilde{\upsilon}_{m_1...m_N}(y)dy$$

Следовательно, применяя метод шагов (см., например [2], [3]), получим:

$$D_{0t}^{\alpha} T_{m_1...m_N}(t) + \mu_{m_1...m_N} a^2 T_{m_1...m_N}(t) = f_{m_1...m_N}(t) - \mu_{m_1...m_N} b^2 (D_{0t}^{l-\alpha}\varphi_l)_{m_1...m_N}(t) \qquad (10)$$

при $0 \leq t \leq \tau$ и начальные условие

$$D_{0t}^{\alpha-i} T_{m_1...m_N}(t)\big|_{t=0} = \varphi_{i,m_1...m_N}, \quad i=1,2,...,l \qquad (11)$$

Решение задачи Коши (10), (11) известно (см., например [6, с. 16-17]) и оно имеет вид

$$_1T_{m_1...m_N}(t) = \sum_{j=1}^{l} \varphi_{j,m_1...m_N} t^{\alpha-j} E_{\alpha,\alpha-j+1}(-\mu_{m_1...m_N} a^2 \cdot t^{\alpha}) +$$

$$+ \int_0^t (t-\xi)^{\alpha-1} \cdot E_{\alpha,\alpha}\left[-\mu_{m_1...m_N} a^2 (t-\xi)^{\alpha}\right] \left[f_{m_1...m_N}(\xi) - \mu_{m_1...m_N} b^2 (D_{0t}^{l-\alpha}\varphi_l)_{m_1...m_N}(\xi)\right]d\xi$$

при $0 \leq t \leq \tau$, где $\qquad (12)$

$$\mu_{m_1...m_N} = \sum_{j=1}^{N} \lambda_{m_j}^2 = \sum_{j=1}^{p}(2m_j + \varphi_j)^2 + \sum_{j=p+1}^{q}(2m_j)^2 + \sum_{j=q+1}^{N}(m_j)^2 \qquad (13)$$

$$E_{\alpha,\beta}(z) = \sum_{q=0}^{\infty} \frac{z^q}{\Gamma(\alpha q + \beta)} \qquad (14)$$

Далее, для отрезка $\tau \leq t \leq 2\tau$ из уравнение (8), получим:

$$D_{0t}^{\alpha} T_{m_1...m_N}(t) + \mu_{m_1...m_N} a^2 T_{m_1...m_N}(t) = f_{m_1...m_N}(t) - \mu_{m_1...m_N} b^2 \cdot {_1T_{m_1...m_N}}(t-\tau). \qquad (15)$$

Из решение (12) получаем начальные условие вида

$$D_{0t}^{\alpha-i} T_{m_1...m_N}(t)\big|_{t=\tau} = {_1T_{i,m_1...m_N}}, \quad i=1,2,...,l \qquad (16)$$

Аналогично, как задачи Коши (10), (11), решение задачи Коши (15), (16) известно и оно имеет вид

$$_2T_{m_1...m_N}(t) = \sum_{j=1}^{l} \varphi_{j,m_1...m_N}(t-\tau)^{\alpha-j} E_{\alpha,\alpha-j+1}(-\mu_{m_1...m_N} a^2 \cdot (t-\tau)^\alpha) +$$

$$+ \int_{\tau}^{t} (t-\xi)^{\alpha-1} \cdot E_{\alpha,\alpha}\left[-\mu_{m_1...m_N} a^2 (t-\xi)^\alpha\right]\left[f_{m_1...m_N}(\xi) - \mu_{m_1...m_N} b^2 {}_1T_{m_1...m_N}(\xi-\tau)\right]d\xi$$

при $\tau \leq t \leq 2\tau$. (17)

Далее, применяя метод шагов для каждого отрезка $n\tau \leq t \leq (n+1)\tau$ из уравнение (8), получим:

$$D_{0t}^\alpha T_{m_1...m_N}(t) + \mu_{m_1...m_N} a^2 T_{m_1...m_N}(t) = f_{m_1...m_N}(t) - \mu_{m_1...m_N} b^2 \cdot {}_nT_{m_1...m_N}(t-n\tau). \quad (18)$$

Соответствующие начальные условие имеет вид

$$D_{0t}^{\alpha-i} T_{m_1...m_N}(t)\big|_{t=n\tau} = {}_nT_{i,m_1...m_N}, \quad i=1,2,...,l \quad (19)$$

Решая задачу Коши (18), (19), получим:

$$_{(n+1)}T_{m_1...m_N}(t) = \sum_{j=1}^{l} \varphi_{j,m_1...m_N}(t-n\tau)^{\alpha-j} E_{\alpha,\alpha-j+1}(-\mu_{m_1...m_N} a^2 \cdot (t-n\tau)^\alpha) +$$

$$+ \int_{n\tau}^{t} (t-\xi)^{\alpha-1} \cdot E_{\alpha,\alpha}\left[-\mu_{m_1...m_N} a^2 (t-\xi)^\alpha\right]\left[f_{m_1...m_N}(\xi) - \mu_{m_1...m_N} b^2 {}_nT_{m_1...m_N}(\xi-n\tau)\right]d\xi$$

при $n\tau \leq t \leq (n+1)\tau$. (20)

Поскольку функции (6) при каждом отрезке $n\tau \leq t \leq (n+1)\tau$ построены в явном виде с помощью формулами (12), (17), (20), то на основании полноты системы собственных функций (7) в $L_2(\Pi)$ нетрудно доказать единственность решения задачи (1)–(3). Пусть $f(x,t) \equiv 0$ и $\varphi_i(x) \equiv 0, i=1,...,l$.

Тогда из формулы (12), (17), (20) и (6) следует, что

$$\int_{\Pi} u(y,t)\tilde{\upsilon}_{m_1...m_N}(y)dy = 0$$

при всех $m_1, m_2, ..., m_N \in \square$ и любом $t \in [0,T]$. Отсюда в силу полноты системы собственных функций (7) в $L_2(\Pi)$ вытекает, что $u(x,t) = 0$ почти всюду в области $\Pi$ при любом $t \in [0,T]$. Как известно, по теореме вложения Соболева функция $u(x,t)$ непрерывна на $\overline{Q}$, то $u(x,t) \equiv 0$ в $\overline{Q}$. Это доказывают единственность решения задачи (1)–(3).

При каждом $t > 0$ функция $u(x,t) \in \overset{0}{W}_2^{s_1,s_2,...,s_N;\theta}(Q)$ по переменной $x$ является функций из класса $\overset{0}{W}_2^{s_1,s_2,...,s_N}(\Pi)$. Поэтому, рассматривая $t>0$ как параметр, решение задачи (1)–(3) будем искать из класса $\overset{0}{W}_2^{s_1,s_2,...,s_N;\theta}(Q)$ в виде суммы ряда по системе собственных функций (7) спектральной задачи (4), (5):

$$u(x,t) = \sum_{m_1=1}^{\infty} \ldots \sum_{m_N=1}^{\infty} T_{m_1\ldots m_N}(t) \cdot \tilde{\upsilon}_{m_1\ldots m_N}(x), \qquad (21)$$

где $\tilde{\upsilon}_{m_1\ldots m_N}(x)$ определяется по формуле (7) и $T_{m_1\ldots m_N}(t)$ определяется по формулами (12), (17) и (20).

После подстановки (12) в (21) мы получим единственное решение задачи (1)–(3) в виде ряда

$$u(x,t) = \sum_{m_1=1}^{\infty} \ldots \sum_{m_N=1}^{\infty} \Bigg[ \sum_{j=1}^{l} \varphi_{j,m_1\ldots m_N} \cdot t^{\alpha-j} E_{\alpha,\alpha-j+1}(-\mu_{m_1\ldots m_N} a^2 \cdot t^\alpha) +$$

$$+ \int_0^t (t-\xi)^{\alpha-1} \cdot E_{\alpha,\alpha}\left[-\mu_{m_1\ldots m_N} a^2 (t-\xi)^\alpha\right]\left[f_{m_1\ldots m_N}(\xi) - \mu_{m_1\ldots m_N} b^2 (D_{0t}^{l-\alpha}\varphi_l)_{m_1\ldots m_N}(\xi)\right]d\xi \times$$

$$\times \tilde{\upsilon}_{m_1\ldots m_N}(x_1,\ldots,x_N)$$
(22)

при $0 \leq \tau \leq t$. Аналогично, после подстановки (17) в (21) мы получим в отрезке $\tau \leq t \leq 2\tau$ единственное решение задачи (1)–(3) в виде ряда

$$u(x,t) = \sum_{m_1=1}^{\infty} \ldots \sum_{m_N=1}^{\infty} \Bigg[ \sum_{j=1}^{l} \varphi_{j,m_1\ldots m_N} \cdot (t-\tau)^{\alpha-j} E_{\alpha,\alpha-j+1}(-\mu_{m_1\ldots m_N} a^2 \cdot (t-\tau)^\alpha) +$$

$$+ \int_\tau^t (t-\xi)^{\alpha-1} \cdot E_{\alpha,\alpha}\left[-\mu_{m_1\ldots m_N} a^2 (t-\xi)^\alpha\right]\left[f_{m_1\ldots m_N}(\xi) - \mu_{m_1\ldots m_N} b^2 \cdot {}_1T_{m_1\ldots m_N}(\xi-\tau)\right]d\xi \times$$

$$\times \tilde{\upsilon}_{m_1\ldots m_N}(x_1,\ldots,x_N).$$
(23)

Следовательно, после подстановки (20) в (21) мы получим в отрезке $n\tau \leq t \leq (n+1)\tau$ единственное решение задачи (1)–(3) в виде ряда

$$u(x,t) = \sum_{m_1=1}^{\infty} \ldots \sum_{m_N=1}^{\infty} \Bigg[ \sum_{j=1}^{l} \varphi_{j,m_1\ldots m_N} \cdot (t-n\tau)^{\alpha-j} E_{\alpha,\alpha-j+1}(-\mu_{m_1\ldots m_N} a^2 \cdot (t-n\tau)^\alpha) +$$

$$+ \int_{n\tau}^t (t-\xi)^{\alpha-1} \cdot E_{\alpha,\alpha}\left[-\mu_{m_1\ldots m_N} a^2 (t-\xi)^\alpha\right]\left[f_{m_1\ldots m_N}(\xi) - \mu_{m_1\ldots m_N} b^2 \cdot {}_nT_{m_1\ldots m_N}(\xi-n\tau)\right]d\xi \times$$

$$\times \tilde{\upsilon}_{m_1\ldots m_N}(x_1,\ldots,x_N).$$
(24)

Поскольку система собственных функций (7) спектральной задачи (4), (5) образует базис Рисса в пространстве Соболева $\overset{0}{H}{}^{s_1,s_2,\ldots,s_N}(\prod)$, то любая функция из этого класса разлагается единственным образом в ряд Фурье сходящийся по норме пространства $H^{s_1,s_2,\ldots,s_N}(\prod)$. Поэтому ряды (22), (23) и (24) сходится в $H^{s_1,s_2,\ldots,s_N}(\prod)$ соответственно в отрезках $t \in [0,\tau]$, $t \in [\tau,2\tau]$ и $t \in [n\tau,(n+1)\tau]$. Выясним условия существования решения из класса $\overset{0}{W}{}_2^{s_1,s_2,...,s_N;\theta}(Q)$. Согласно теореме 1, система собственных функций (7) спектральной задачи (4), (5) образует базис Рисса в пространствах Соболева $\overset{0}{H}{}^{s_1,s_2,...,s_N}(\Pi)$ и $\overset{0}{W}{}^{s_1,s_2,...,s_N}(\Pi)$. Поэтому справедливо неравенство

$$\|u(x,t)\|^2_{C^{2,2,...,2}(\Pi)} \le c_8 \|u(x,t)\|^2_{H^{s_1,s_2,...,s_N}(\Pi)} \le$$

$$\le c_9 \sum_{m_1=1}^{\infty}...\sum_{m_N=1}^{\infty} \left| \sum_{j=1}^{l} \varphi_{j,m_1...m_N} t^{\alpha-j} E_{\alpha,\alpha-j+1}(-\mu_{m_1...m_N} a^2 \cdot t^{\alpha}) + \right.$$

$$\left. + \int_0^t (t-\xi)^{\alpha-1} E_{\alpha,\alpha}\left[-\mu_{m_1...m_N} a^2(t-\xi)^{\alpha}\right]\left[f_{m_1...m_N}(\xi) - \mu_{m_1...m_N} b^2 (D_{0t}^{l-\alpha}\varphi_l)_{m_1...m_N}(\xi)\right]d\xi \right|^2 \times$$

$$\times \prod_{k=1}^{N}(1+\lambda_{m_k}^{2s_k}) < \infty \qquad (25)$$

при $0 \le t \le \tau$. Далее,

$$\|u(x,t)\|^2_{C^{2,2,...,2}(\Pi)} \le c_8 \|u(x,t)\|^2_{H^{s_1,s_2,...,s_N}(\Pi)} \le$$

$$\le c_9 \sum_{m_1=1}^{\infty}...\sum_{m_N=1}^{\infty} \left| \sum_{j=1}^{l} \varphi_{j,m_1...m_N} (t-\tau)^{\alpha-j} E_{\alpha,\alpha-j+1}(-\mu_{m_1...m_N} a^2 \cdot (t-\tau)^{\alpha}) + \right.$$

$$\left. + \int_\tau^t (t-\xi)^{\alpha-1} E_{\alpha,\alpha}\left[-\mu_{m_1...m_N} a^2(t-\xi)^{\alpha}\right]\left[f_{m_1...m_N}(\xi) - \mu_{m_1...m_N} b^2 \,_1T_{m_1...m_N}(\xi-\tau)\right]d\xi \right|^2 \times$$

$$\times \prod_{k=1}^{N}(1+\lambda_{m_k}^{2s_k}) < \infty, \qquad (26)$$

при $\tau \le t \le 2\tau$. Следовательно,

$$\|u(x,t)\|^2_{C^{2,2,...,2}(\Pi)} \le c_8 \|u(x,t)\|^2_{H^{s_1,s_2,...,s_N}(\Pi)} \le$$

$$\le c_9 \sum_{m_1=1}^{\infty}...\sum_{m_N=1}^{\infty} \left| \sum_{j=1}^{l} \varphi_{j,m_1...m_N} (t-n\tau)^{\alpha-j} E_{\alpha,\alpha-j+1}(-\mu_{m_1...m_N} a^2 \cdot (t-n\tau)^{\alpha}) + \right.$$

$$\left. + \int_{n\tau}^t (t-\xi)^{\alpha-1} E_{\alpha,\alpha}\left[-\mu_{m_1...m_N} a^2(t-\xi)^{\alpha}\right]\left[f_{m_1...m_N}(\xi) - \mu_{m_1...m_N} b^2 \,_nT_{m_1...m_N}(\xi-n\tau)\right]d\xi \right|^2 \times$$

$$\times \prod_{k=1}^{N}(1+\lambda_{m_k}^{2s_k}) < \infty \qquad (27)$$

при $n\tau \le t \le (n+1)\tau$.

По теореме вложения Соболева условие (25), (26) и (27) является достаточным условием существования регулярного решения задачи (1)–(3) из класса $\overset{0}{W}{}_2^{s_1,s_2,...,s_N;\theta}(Q)$ с показателем $s_1 = s_2 = ... = s_N = 2 + \dfrac{N}{2}, \theta = -[-\alpha]$.

Если $0 < \alpha < 2$, то с учетом оценки функции Миттаг–Леффлера (см., например [7, с. 35])

$$\left|E_{\alpha,\alpha-j+1}(\mu_{m_1...m_N}t^\alpha)\right| \leq \frac{C}{1+\left|\mu_{m_1...m_N}t^\alpha\right|}, \quad \left|E_{\alpha,\alpha}(\mu_{m_1...m_N}(t-\tau)^\alpha)\right| \leq \frac{C}{1+\left|\mu_{m_1...m_N}(t-\tau)^\alpha\right|}$$

можно упростить достаточные условия (25), (26) и (27) существования регулярного решения задачи (1)–(3) из класса $\overset{0}{W}{}_2^{s_1,s_2,...,s_N;\theta}(Q)$ с показателем $s_1 = s_2 = ... = s_N = 2 + \dfrac{N}{2}, \theta = -[-\alpha]$.

Если для любого $j = 1, 2, ..., l$ и $0 \leq t \leq \tau$

$$\sum_{m_1=1}^{\infty} ... \sum_{m_N=1}^{\infty} \left|\varphi_{j,m_1...m_N}\right|^2 \prod_{k=1}^{N}(1+\lambda_{m_k}^{s_k}) < \infty,$$

$$\sum_{m_1=1}^{\infty} ... \sum_{m_N=1}^{\infty} \left|\int_0^t (t-\xi)^{\alpha-1}\left|f_{m_1...m_N}(\xi) - \mu_{m_1...m_N}b^2 (D_{0t}^{l-\alpha}\varphi_l)_{m_1...m_N}(\xi)\right|d\xi\right|^2 \cdot \prod_{k=1}^{N}(1+\lambda_{m_k}^{s_k}) < \infty \quad (28)$$

то условия (25) будет выполнено.

Далее, если для любого $j = 1, 2, ..., l$ и $\tau \leq t \leq 2\tau$

$$\sum_{m_1=1}^{\infty} ... \sum_{m_N=1}^{\infty} \left|\varphi_{j,m_1...m_N}\right|^2 \prod_{k=1}^{N}(1+\lambda_{m_k}^{s_k}) < \infty,$$

$$\sum_{m_1=1}^{\infty} ... \sum_{m_N=1}^{\infty} \left|\int_0^t (t-\xi)^{\alpha-1}\left|f_{m_1...m_N}(\xi) - \mu_{m_1...m_N}b^2 \left(D_{0t}^{l-\alpha}\varphi_l\right)_{m_1...m_N}(\xi)\right|d\xi\right|^2 \cdot \prod_{k=1}^{N}\left(1+\lambda_{m_k}^{2s_k}\right) < \infty,$$

$$\sum_{m_1=1}^{\infty} ... \sum_{m_N=1}^{\infty} \left|\int_\tau^t (t-\xi)^{\alpha-1}\left|f_{m_1...m_N}(\xi) - \mu_{m_1...m_N}b^2 \cdot {}_1T_{m_1...m_N}(\xi-\tau)\right|d\xi\right|^2 \cdot \prod_{k=1}^{N}\left(1+\lambda_{m_k}^{2s_k}\right) < \infty,$$

(29)

то условие (26) будет выполнено при $\tau \leq t \leq 2\tau$.

Следовательно, если для любого $j = 1, 2, ..., l$ и $n\tau \leq t \leq (n+1)\tau$

$$\sum_{m_1=1}^{\infty} ... \sum_{m_N=1}^{\infty} \left|\varphi_{j,m_1...m_N}\right|^2 \prod_{k=1}^{N}(1+\lambda_{m_k}^{s_k}) < \infty,$$

$$\sum_{m_1=1}^{\infty} ... \sum_{m_N=1}^{\infty} \left|\int_0^t (t-\xi)^{\alpha-1}\left|f_{m_1...m_N}(\xi) - \mu_{m_1...m_N}b^2 \left(D_{0t}^{l-\alpha}\varphi_l\right)_{m_1...m_N}(\xi)\right|d\xi\right|^2 \cdot \prod_{k=1}^{N}\left(1+\lambda_{m_k}^{2s_k}\right) < \infty,$$

$$\sum_{m_1=1}^{\infty} ... \sum_{m_N=1}^{\infty} \left|\int_\tau^t (t-\xi)^{\alpha-1}\left|f_{m_1...m_N}(\xi) - \mu_{m_1...m_N}b^2 \cdot {}_1T_{m_1...m_N}(\xi-\tau)\right|d\xi\right|^2 \cdot \prod_{k=1}^{N}\left(1+\lambda_{m_k}^{2s_k}\right) < \infty,$$

................................................,

$$\sum_{m_1=1}^{\infty} ... \sum_{m_N=1}^{\infty} \left|\int_{n\tau}^t (t-\xi)^{\alpha-1}\left|f_{m_1...m_N}(\xi) - \mu_{m_1...m_N}b^2 \cdot {}_nT_{m_1...m_N}(\xi-n\tau)\right|d\xi\right|^2 \cdot \prod_{k=1}^{N}\left(1+\lambda_{m_k}^{2s_k}\right) < \infty,$$

(30)

то условие (27) будет выполнено.

Справедлива следующая

**Теорема 5.** *Пусть начальные функций $\varphi_i(x)$, $i=1,2,...,l$, и правую часть $f(x,t)$ удовлетворяет условию (30) при каждом $n\tau \leq t \leq (n+1)\tau$. $n=0,1,...,\left[\dfrac{T}{\tau}\right]$. Тогда регулярные решение задачи (1)–(3) из класса $\overset{0}{W}{}_2^{s_1,s_2,...,s_N;\theta}(Q)$ с показателем $s_1=s_2=...=s_N=2+\dfrac{N}{2}$, $\theta=-[-\alpha]$ существует, единственно и представляется в виде ряда (22), (23) и (24), где коэффициенты определяется по формулами (12), (17) и (20).*

## ЛИТЕРАТУРА